\documentstyle[11pt]{article}
\input{amssym.def}
\input epsf

\newtheorem{thm}{Theorem}[section]
\newtheorem{lemma}[thm]{Lemma}
\newtheorem{prop}[thm]{Proposition}

\newtheorem{cor}[thm]{Corollary}

\newtheorem{ex}[thm]{Example}

\newenvironment{pf}{\bf Proof:\rm}{\begin{flushright}$\Box$\end{flushright}}
\newcommand{\Qed}{\begin{flushright}$\Box$\end{flushright}}

\newcommand{\as}{{\rm{as}}}

\def\rank{\mathop{\rm rank}\nolimits}
\def\ker{\mathop{\rm Ker}\nolimits}

\newcommand{\C}{{\Bbb C}}

\newcommand{\R}{{\Bbb R}}

\newcommand{\HH}{{\Bbb H}}
\newcommand{\cS}{{\cal S}}
\newcommand{\cA}{{\cal A}}

\newcommand{\cK}{{\cal K}}
\newcommand{\cL}{{\cal L}}
\newcommand{\cN}{{\cal N}}

\newcommand{\cP}{{\cal P}}

\newcommand{\cX}{{\cal X}}

\newcommand{\cV}{{\cal V}}
\newcommand{\cW}{{\cal W}}

\newcommand{\pF}{{\cal F}}

\newcommand{\fg}{{\frak g}}

\newcommand{\fk}{{\frak k}}
\newcommand{\fp}{{\frak p}}
\newcommand{\bb}{{\backslash\!\backslash}}

\begin{document}
\bibliographystyle{plain}

\title{On the Specialization to the Asymptotic Cone}
\author{Mikhail Grinberg}
\maketitle

\section{Introduction}

Let $V \cong \C^d$ be a complex vector space, and $X \subset V$
a connected, smooth, closed subvariety.  Consider the family
of scalings $\{ \lambda X \}_{\lambda \in \C^*}$.  The limit
$\as (X) := \lim_{\lambda \to 0} \; \lambda X$ is a conic
subvariety of $V$, called the {\em asymptotic cone} of $X$.  The
asymptotic cone is equipped with a nearby cycles perverse sheaf
$P = P^{}_X$, coming from the specialization of $X$ to $\as (X)$
(we give precise definitions in Section 2.2).  The sheaf $P$
is $\C^*$-conic.

Given a Hermitian inner product on $V$, we may consider, for any
$x \in X$, the angle $\angle (x, T_x X) \in [0, \pi/2]$ between the
vector $x \in V$ and the subspace $T_x X \subset T_x V \cong V$.
We say that $X$ is {\em transverse to infinity} if for some
(equivalently for any) inner product on $V$, there exists
a constant $k > 0$, such that for any $x \in X$, we have:
$$\angle (x, T_x X) < \frac{k}{|\!| x |\!|} \, .$$
For example, a hyperbola in the plane is transverse to infinity,
but a parabola is not.  Our main result is the following.

\begin{thm}\label{main}
Assume $X \subset V$ are as above, and $X$ is transverse to
infinity.  Let $T_X^* V \subset T^* V \cong V \times V^*$ be the
conormal bundle to $X$, and $p_2 : V \times V^* \to V^*$ be the
projection on the second factor.  Let $Y = \overline{p_2 (T_X^* V)}$;
it is an irreducible cone in $V^*$.  Then we have:
$$\pF \, P^{}_X \cong {\cal IC} (Y^\circ, \cL),$$
where $\pF$ is the geometric Fourier transform functor, $\cL$ is
a local system on some Zariski open subset $Y^\circ \subset Y$, and
the right hand side is the intersection homology extension of $\cL$. 
\end{thm}

The reader will, at a glance, place this result in the context
of the classical Gauss map / projective duality ideas.  In the
modern literature (known to the author) it is perhaps closest
to Jean-Luc Brylinski's paper [Br], which develops the Lefschetz
pencil and projective duality themes in the context of constructible
sheaves.  Here are some situations where the theorem applies.

\begin{ex}\label{generic_homogeneous}{\em
Let $f : V \to \C$ be a generic homogeneous polynomial of degree
$m$, and $X = f^{-1} (1)$.  Then $\as (X) = f^{-1} (0)$, and $P$
is the sheaf of nearby cycles of $f$.  In this case, $Y = V^*$
and $\rank \cL = m (m-1)^{d-1}$.}
\end{ex}

\begin{ex}\label{linear_forms}{\em
Let $f  = \prod l_i: V \to \C$ be a product of (finitely many)
linear forms $l_i \in V^*$, and $X = f^{-1} (1)$.  As in the
previous example, $\as (X) = f^{-1} (0)$, and $P$ is the sheaf of
nearby cycles of $f$.  In this case, $Y = (\cap \ker l_i)^\perp$.
We will discuss this example in more detail in Section 6.  In
particular, we will give an algorithm for computing the rank of
the local system $\cL$.}
\end{ex}

\begin{ex}\label{springer}{\em
Let $V = \fg$ be a complex semisimple Lie algebra, and
$X \subset \fg$ be a semisimple (i.e., closed) adjoint orbit.
Then $X$ is transverse to infinity, and $\as (X) \subset \fg$
is a closed union of nilpotent orbits.  If $X$ is a regular
orbit, then $\as (X) = \cN$, the full nilcone in $\fg$, and $P$
is the famous Springer perverse sheaf, also known as the nearby
cycles of the adjoint quotient map $f : \fg \to G \bb \fg$
(see [Sp], [Sl], [BM], [M]).  The claim of Theorem \ref{main} in
this case is due independently to Hotta-Kashiwara [HK] and
Ginzburg [Gi] (see also [Br]).  In a way, Theorem \ref{main} is
saying that the Fourier transform description of Springer theory
is independent of the existence of the simultaneous resolution
for the fibers of $f$ (the celebrated Grothendieck-Springer
resolution).}
\end{ex}

\begin{ex}\label{sym_sp}{\em
A natural generalization of Example \ref{springer} is the
following.  Let $\fg = \fk \oplus \fp$ be a Cartan decomposition
of a complex semisimple Lie algebra.  Take $V = \fp$, and let
$X \subset \fp$ be a closed $K$-orbit.  Then $X$ is transverse to
infinity.  If $X$ is a regular orbit, then $\as (X)$ is the
nilcone in $\fp$, and Theorem \ref{main} can be used to give a
generalization of Springer theory to symmetric spaces (where no
simultaneous resolution is available).}
\end{ex}

A further generalization of Examples \ref{springer} and
\ref{sym_sp} is given by the so-called polar representations
of Dadok and Kac (see [DK]).  Their geometry forms the subject of
the paper [Gr1].

We should note that the local system $\cL$ appearing in Theorem
\ref{main} typically will not be semisimple.  Thus, even in the
cases where $\cL$ can be explicitly computed (e.g., Example
\ref{sym_sp}), Theorem \ref{main} may fall short of giving the
complete structure of $P$ as a perverse sheaf.  This is because
intersection homology is not an exact functor from local systems
to perverse sheaves.  However, the theorem is saying that, in some
sense, the whole nearby cycles sheaf is encoded in a single
local system.  In particular, the endomorphism rings of $P$ and
$\cL$ are the same.

The paper is organized as follows.  We fix the notations in
Section 2.  Sections 3-5 are devoted to the proof of Theorem
\ref{main}.  Lastly, Section 6 is concerned with products of
linear forms (Example \ref{linear_forms} above).  The reader
interested in symmetric spaces or, more generally, polar
representations is referred to [Gr1] and [Gr2].

I would like to thank the IAS, Princeton, where most of this work
was done, for its hospitality, and the Fannie and John Hertz
foundation for financial support.  Discussions with Sam Evens,
Victor Ginzburg, Ian Grojnowski, David Massey, and Ivan Mirkovi\'c
have been of great value to me.  I am also grateful to Mark Goresky,
David Kazhdan, and Wilfried Schmid for their interest and
encouragement.  Numerous conversations with Tom Braden have been a
constant source of help and motivation.  Finally, I would like to
thank Robert MacPherson for introducing me to nearby cycles, and for
his guidance and support throughout this work.

\section{Notations}

\subsection{Sheaves and Functors}
We will say {\em sheaf} to mean {\em complex of sheaves}
throughout; all our sheaves will be sheaves of $\C$ vector
spaces.  Given a map $g : X \to Y$, the symbols $g_*$, $g_!$
will always denote the {\em derived} push-forward functors.
All perverse sheaves and intersection homology will be taken
with respect to the middle perversity (see [GM1], [BBD]);
we use the shift conventions of [BBD].  Given a sheaf $A$ on
$X$, and a pair of closed subspaces $Z \subset Y \subset X$,
we will write $\HH^k(Y, Z; A)$ for the hypercohomology group
$\HH^k(j_! \,i^* A)$, where $i : Y \setminus Z \to X$ and
$j : Y \setminus Z \to Y$ are the inclusion maps.  We call
$\HH^k(Y, Z; A)$ the relative hypercohomology of $A$.  For
an analytic function $f : M \to \C$, we use the notation of
[KS, Chapter 8.6] for the nearby and vanishing cycles
functors $\psi^{}_f$, $\phi^{}_f$.

When $V$ is a $\C$ vector space, we denote by $\cP^{}_{\C^*}
\, (V)$ the category of $\C^*$-conic perverse sheaves on $V$,
and by $\pF : \cP^{}_{\C^*} \, (V) \to \cP^{}_{\C^*} \, (V^*)$
the (shifted) Fourier transform functor.  In the notation
of [KS, Chapter 3.7], we have $\pF \, P = P \, \hat{} \;
[\dim \, V]$.  To avoid cumbersome notation, we use the
following shorthand: if $P$ is a conic perverse sheaf on a
closed conic subvariety $X \subset V$, and $j^{}_X : X \to V$
is the inclusion, then we write $\pF \, P$ instead of
$\pF \circ (j^{}_X)_* \, P$.

\subsection{The Asymptotic Cone}

Let $V \cong \C^d$ be a complex vector space, and $X \subset V$
a connected, smooth, closed subvariety, as in Section 1.
We denote by $\bar V$ the standard projective compactification
of $V$, and by $\bar X$ the closure of $X$ in $\bar V$.  Set
$V^\infty = \bar V \setminus V$, and $X^\infty = \bar X \cap
V^\infty$.  The asymptotic cone $\as (X) \subset V$ is defined
as the affine cone over $X^\infty$.

Another way to define $\as (X)$ is as follows.  Let 
$\tilde X^\circ = \{ (\lambda, \tilde x) \in \C^* \times V \, |
\, \tilde x \in \lambda X \}$, and $\tilde X$ be the closure
of $\tilde X^\circ$ in $\C \times V$.  Write $\tilde f :
\tilde X \to \C$ for the projection on the first factor.
Then $\as (X) = \tilde f^{-1} (0)$.  The nearby cycles sheaf
$P = P^{}_X$ appearing in Theorem \ref{main} is defined by
$P = \psi_{\tilde f} \, \C_{\tilde X} [n]$, where $n = \dim X$.
Here and in the rest of the paper the symbol $\dim$ denotes
the {\em complex} dimension.

\section{Outline of Proof of Theorem \ref{main}}

The main part of the proof of Theorem \ref{main} is deformation
arguments, \`a la ``moving the wall'' in [GM3].  We give an
outline of the proof in this section, postponing the deformation
arguments to Section 4, and deferring one other technical argument
to Section 5.

\subsection{The Stalks of the Fourier Transform}

Our proof of Theorem \ref{main} will be based on the original
definition of the intersection homology complex by stalk and
costalk cohomology vanishing conditions  [GM1, Theorem 4.1].  We
therefore begin with a lemma that identifies the stalks of the
Fourier transform.

\begin{lemma}\label{ftstalk}
Let $V$ be a Hermitian complex vector space of dimension $d$,
$Q \in \cP^{}_{\C^*} (V)$ an algebraically constructible conic
perverse sheaf, and $l \in V^*$ a covector.  For $v \in V$, set
$\xi (v) = \mbox{\em Re} \, (l(v))$.  Fix positive numbers
$0 < \xi^{}_0 \ll \eta^{}_0$.  Then the stalk cohomology
$H^i_l \, (\pF \, Q)$ is given by a natural isomorphism:
$$H^i_l \, (\pF \, Q) \cong \HH^{i+d} \, ( \{ \, v \in V \, | \,
\xi (v) \leq \xi^{}_0 \, \}, \, \{ \, |\!| v |\!| \geq \eta^{}_0
\, \}; \, Q ),$$
where the right hand side is a relative hypercohomology group
of $Q$.
\end{lemma}
\begin{pf}
This follows from the definition of the Fourier transform
[KS, Chapter 3.7].
\end{pf}

Let now $V \supset X$ be as in Theorem \ref{main}.  Fix a
Hermitian metric on $V$.  Set $n = \dim \, X$.

\begin{cor}\label{ftstalknc}
Fix $l \in V^*$, and let $\xi = \mbox{\em Re} \, (l)$.  Also fix
positive numbers $0 < \lambda \ll \xi^{}_0 \ll \eta^{}_0$, chosen
in decreasing order.  Let $\lambda X$ be the translate of $X$ by
the scalar $\lambda$.  Then we have:
$$H^i_l \, (\pF \, P) \cong H^{i+d+n} \, ( \{ \, \tilde x \in
\lambda X \, | \, \xi (\tilde x) \leq \xi^{}_0 \, \}, \, \{ \,
|\!| \tilde x |\!| \geq \eta^{}_0 \, \}; \, \C ),$$
where the right hand side is an ordinary relative cohomology
group of $\lambda X$.
\end{cor}
\begin{pf}
This follows from Lemma \ref{ftstalk} and the definition of nearby
cycles.
\end{pf}

Note that in Corollary \ref{ftstalknc} the number $\eta^{}_0$ is
chosen first, then $\xi^{}_0$ is chosen small relative to
$\eta^{}_0$, then $\lambda$ is chosen small relative to
$\xi^{}_0$.  Note also that the group in the right hand side is
invariant under the simultaneous multiplication of the numbers
$\lambda, \xi^{}_0, \eta^{}_0$ by a positive scalar. 

\begin{prop}\label{odlim}
The right hand side of the isomorphism of Corollary
\ref{ftstalknc} will not change if one chooses the numbers
$\lambda, \xi^{}_0, \eta^{}_0$ in the opposite order.  Taking
$\lambda = 1$, we have:
$$H^i_l \, (\pF \, P) \cong H^{i+d+n} \, ( \{ \, x \in X \, | \,
\xi (x) \leq \xi^{}_0 \, \}, \, \{ \, |\!| x |\!| \geq \eta^{}_0
\, \}; \, \C ),$$
where the numbers $1 \ll \xi^{}_0 \ll \eta^{}_0$ are chosen large
in increasing order.
\end{prop}

We will give a proof of Proposition \ref{odlim} in Section 4.2.
It is the key step of the proof of Theorem \ref{main}, and the
main place where the transverse-to-infinity assumption is used.

\subsection{A Decomposition of $V^*$}

To verify the definition of intersection homology for the sheaf
$\pF \, P$, we are going to fix a decomposition of $V^*$ into
algebraic manifolds.  For this, we first fix a stratification
$\cS$ of $\as (X)$.  Let $\as (X) = \bigcup_{S \in \cS} S$ be an
algebraic Whitney stratification satisfying the following three
conditions.

(i)   $\cS$ is conic, i.e., each $S \in \cS$ is $\C^*$-invariant.

(ii)  Thom's $A_{\tilde f}$ condition holds for the pair
$(\tilde X^\circ, S)$, for each $S \in \cS$.

(iii) Let $\cS^\circ = \cS \setminus \{\{ 0 \}\}$.  For
$S \in \cS^\circ$, let $S^\infty \subset X^\infty$ be the
projectivization of $S$.  Then the decomposition $\bar X =
X \cup \bigcup_{S \in \cS^\circ} S^\infty$ is a Whitney
stratification.

The reader is referred to [Hi] and [GM3, Part I, Chapter 1]
for a discussion of the $A_{\tilde f}$ condition and of the
Whitney conditions.  The existence of a stratification $\cS$
satisfying (i)-(iii) above follows from the general results
contained in these references.  

Let $\Lambda_X$ denote the conormal bundle $T_X^* V \subset V
\times V^*$, and $\Lambda_{\as (X)} \subset V \times V^*$
be the conormal variety to the stratification $\cS$.
Recall the projection $p_2 : V \times V^* \to V^*$.
We now fix a finite decomposition $V^* = \bigcup_{W \in \cW} W$
of $V^*$ into connected algebraic manifolds, such that for
each $W \in \cW$, the dimensions
$$\dim \, p_2^{-1} (l) \cap \Lambda_X \;\; \mbox{and} \;\;
\dim \, p_2^{-1} (l) \cap \Lambda_{\as (X)}$$
are independent of $l$, for $l \in W$.  Note that one of the
pieces $W \in \cW$ must be an open subset of $Y = \overline
{p_2 (\Lambda_X)}$.  We denote this piece by $W_0$.  The set
$Y^\circ$ appearing in Theorem \ref{main} will be a subset of
$W_0$.  Note that it does not matter for us whether $\cW$ is a
stratification.  For $W \in \cW$, put $d (W) = \dim \, W$.

\begin{prop}\label{deligneax}
(i)  If $l \in W \neq W_0$, then $H^i_l \, (\pF \, P) = 0$, for
$i \geq - d (W)$.

(ii) If $l \in W_0$, then $H^i_l \, (\pF \, P) = 0$, for
$i > - d (W_0)$.
\end{prop}

Sections 3.3 and 3.4 will be devoted to the proof of Proposition
\ref{deligneax}.  Theorem \ref{main} is an immediate
consequence.  Indeed, Proposition \ref{deligneax} is just a
rewriting of half of the axioms for the intersection homology sheaf
[GM1, Theorem 4.1].  The other half, dealing with costalk cohomology,
follows form the fact that $P$, and therefore $\pF \, P$, is Verdier
self-dual.

\subsection{A Partial Compactification of $X$}

Our proof of Proposition \ref{deligneax} will be based on the study
of Morse theory of the function $l$ on $X$.  The difficulty
in analyzing this Morse theory is that $l |_X$ is not proper.  For
this reason, we introduce a partial compactification $\hat X =
\hat X_l$ of $X$, depending on $l \in V^*$.  In the following
construction $l \in V^*$ is fixed, and we assume $l \neq 0$.
(Case $l = 0$ of Proposition \ref{deligneax} is trivial.)

Let $\Delta \subset V$ be the kernel of $l$, and $L \subset V$ be
any line complementary to $\Delta$.  We have $V = \Delta \oplus L$.
Take the standard projective compactification $\bar \Delta$ of
$\Delta$, and let $\hat V  = \hat V_l = \bar \Delta \times L$.  It
is not hard to check that the space $\hat V$ is canonically
independent of the choice of the line $L$.  Note that $l : V \to \C$
extends to a proper algebraic function $\hat l : \hat V \to \C$.
Let $\hat X = \hat X_l$ be the closure of $X$ in $\hat V$.

Put $\hat V^\infty = \hat V \setminus V$, and $\hat X^\infty =
\hat X \setminus X$.  Any point of $\hat V^\infty$ can be written
as a pair $(\Gamma, z)$, where $\Gamma$ is a one dimensional
subspace of $\Delta$, and $z = \hat l (\Gamma, z) \in \C$.  Note
that if $(\Gamma, z) \in \hat X^\infty$, then $\Gamma \subset
\as (X)$.  For each $S \in \cS^\circ$, let $\hat S^{\infty, \circ}
\subset \hat X^{\infty}$ be the set of all pairs $(\Gamma, z)$ such
that $\Gamma \setminus 0 \subset S$, and the restriction $l \, |_S$
is non-critical along $\Gamma \setminus 0$; it is a manifold of
dimension $\dim \, S - 1$.  The following is a straightforward
verification.

\begin{lemma}\label{noncrit}
Let $(\Gamma, z) \in \hat S^{\infty, \circ}$, then
$(\Gamma, z') \in \hat S^{\infty, \circ}$ for any $z' \in \C$.
\end{lemma}
\Qed

Write $\hat X^{\infty, \circ} = \bigcup_{S \in \cS^\circ}
\hat S^{\infty, \circ}$, and $\hat X^\circ = X \cup
\hat X^{\infty, \circ}$.  Note that $\hat X^\circ$ is Zariski
open in $\hat X$.  For $S \in \cS^\circ$, let $\cK (S)$
be the set of connected components of $\hat S^{\infty, \circ}$.

\begin{prop}\label{whitney}
The decomposition
$$\hat X^\circ = X \, \cup \; \bigcup_{S \in \cS^\circ \atop
K \in \cK (S)} K$$
is a Whitney stratification.
\end{prop} 

We defer the proof of Proposition \ref{whitney} to Section 5.
The significance of the compactification $\hat X$ for our
problem is revealed by the following proposition.  Let
$j : X \to \hat X$ be the inclusion map.

\begin{prop}\label{auxmtw}
The statement of Proposition \ref{odlim} can be further modified
as follows:
$$H^i_l \, (\pF \, P) \cong \HH^{i+d} \, ( \{ \, x \in \hat X
\, | \, \hat \xi (x) \leq \xi^{}_0 \, \}, \, \{ \, |\!| \hat l
(x) |\!| \geq 2 \xi^{}_0 \, \}; \; j^{}_! \, \C^{}_X \, [n] ),$$
where $\xi^{}_0 \gg 1$, and $\hat \xi (x) = \mbox{\em Re} \,
(\hat l(x))$.
\end{prop}

We defer the proof of this to Section 4.3.

\subsection{Proof of Proposition \ref{deligneax}}

We are now prepared to give a proof of Proposition
\ref{deligneax}, completing the outline of the proof of
Theorem \ref{main}.  Our argument is based on the following
general result.

\begin{lemma}\label{dimbound}
Let $A$ be a complex algebraic variety, $\cA$ an algebraic
Whitney stratification of $A$, and $Q$ a perverse sheaf
on $A$, constructible with respect to $\cA$.  Let
$g : A \to \C$ be a proper algebraic function, $z^{}_0$ a
stratified critical value of $g$, and $Z^{}_0 \subset g^{-1}
(z^{}_0)$ the set of stratified critical points of $g$ above
$z^{}_0$.  Fix a small number $0 < \epsilon \ll 1$.  Then the
following bound holds:
$$\HH^i \, ( \{ \, a \in A \; | \; {\em Re} \,
(g (a) - z^{}_0) \leq \epsilon \, \}, \, \{ \,
| g(a) - z^{}_0 | \geq 2 \epsilon \, \}; \; Q) = 0,$$
for $i > \dim \, Z^{}_0$.
\end{lemma}
\begin{pf}
Let $\phi_{g - z_0} \, Q$ be the vanishing cycles of $Q$ with
respect to the function $g - z_0$; it is a perverse sheaf on
$Z_0$ [KS, Corollary 10.3.13].  The hypercohomology group in
the statement of the lemma is the dual of
$\HH^{-i} (Z_0, \, \phi_{g - z_0} \, Q)$.  The lemma follows.
\end{pf}

We are going to apply Lemma \ref{dimbound} to the situation
where $A = \hat X$ and $g = \hat l \, | _{\hat X}$.
Fix an algebraic Whitney stratification $\hat \cX$ of $\hat X$,
such that each stratum of the stratification of $\hat X^\circ$
described in Proposition \ref{whitney} is also a stratum of
$\hat \cX$ (see [GM3, Part I, Theorem 1.7] for a proof that
such an $\hat \cX$ exists).

\begin{lemma}\label{extacdiv}
The extension by zero $j^{}_! \, \C^{}_X \, [n]$,
appearing in Proposition \ref{auxmtw}, is a perverse sheaf on
$\hat X$, constructible with respect to $\hat \cX$.
\end{lemma}
\begin{pf}
By [BBD, Corollary 4.1.10], extending a perverse sheaf by zero
across a divisor always gives a perverse sheaf.
Constructibility with respect to $\hat \cX$ follows from the
local topological triviality of Whitney stratifications 
[GM3, Part I, Chapter 1.4].
\end{pf}

Let $Z \subset \hat X$ be the stratified critical points of
$\hat l \, | _{\hat X}$ with respect to $\hat \cX$.
Recall the decomposition $V^* = \bigcup_{W \in \cW} W$ of
Section 3.2.  Let $W$ be the piece containing $l$.

\begin{lemma}\label{bounddim}
(i)  If $W \neq W_0$, then $\dim \, Z < d - d (W)$.

(ii) If $W = W_0$, then $\dim \, Z = d - d (W_0)$.
\end{lemma}
\begin{pf}
We will consider the sets $Z \cap X$ and $Z \cap \hat X^\infty$
separately.  Look first at $Z \cap X$.  By definition,
$Z \cap X \cong p_2^{-1} (l) \cap \Lambda_X$.  By assumption,
the dimension of $p_2^{-1} (l') \cap \Lambda_X$ is the same for
all $l' \in W$.  This implies:
$$\dim \, p_2^{-1} (W) \cap \Lambda_X =
\dim \, p_2^{-1} (l) \cap \Lambda_X + d(W).$$
Note now that $p_2^{-1} (W_0) \cap \Lambda_X$ is Zariski open
in $\Lambda_X$, which is irreducible of dimension $d$.  The
required bounds on the dimension of $Z \cap X$ follow.

Consider now the set $Z \cap \hat X^\infty$.  We will show that 
$$\dim \, Z \cap \hat X^\infty < d - d (W),$$ 
in both cases of the lemma.  The image $\hat l (Z) \subset \C$ is
a finite set of points, so it will suffice to show that
\begin{equation}\label{DD}
\dim \, Z \cap \hat X^\infty \cap \hat l^{-1} (z) < d - d (W),
\end{equation}
for every $z \in \C$.  By Lemma \ref{noncrit}, we have
$Z \cap \hat X^\infty \subset \hat X^\infty \setminus
\hat X^{\infty, \circ}$.  Let $(\Gamma, z) \in \hat X^\infty
\setminus \hat X^{\infty, \circ}$, where $\Gamma$ is a line in
$\as (X) \cap \Delta$.  Then, for any $v \in \Gamma \setminus 0$,
we have $(v, l)  \in  \Lambda_{\as (X)} \subset V \times V^*$.
It follows that:
$$\dim \, Z \cap \hat X^\infty \cap \hat l^{-1} (z)
\leq \dim \, p_2^{-1} (l) \cap \Lambda_{\as (X)} - 1.$$ 
But the dimension of $p_2^{-1} (l') \cap \Lambda_{\as (X)}$
is the same for all $l' \in W$.  Therefore,
$$\dim \, p_2^{-1} (l) \cap \Lambda_{\as (X)} + d (W) \leq
\dim \, \Lambda_{\as (X)} = d.$$
Inequality (\ref{DD}) follows.
\end{pf}

Proposition \ref{deligneax} now follows from Proposition
\ref{auxmtw} and Lemmas \ref{dimbound} - \ref{bounddim}. 
More precisely, we have to iterate Lemma \ref{dimbound}
for every critical value $z$ of $\hat l \, |_{\hat X}$.

\section{Moving the Wall}

The purpose of this section is to give proofs of Propositions
\ref{odlim} and \ref{auxmtw}.  Both are accomplished by moving the
wall arguments (see [GM3, Part I, Chapter 4] for a discussion of
the moving the wall technique), and both use the assumption that
$X$ is transverse to infinity in a crucial way.

\subsection{Preliminaries}

We continue with a fixed non-zero $l \in V^*$.  Pick a Hermitian
metric on $V$, such that $|\!| l |\!| = 1$.  We introduce
several constants reflecting the geometry of $X$, and the choice
of the covector $l$.  First, there is the number $k > 0$ appearing
in the transverse-to-infinity condition:
\begin{equation}\label{tti}
\angle (x, T_x X) < \frac{k}{|\!| x |\!|},
\end{equation}
for all $x \in X$.  We will use a shorthand $c^{}_1 = 1000 \, k$.
Recall that we denote by $Z \subset \hat X$ the stratified
critical locus of $\hat l \, |_{\hat X}$ with respect to the
stratification $\hat \cX$. The set $\hat l (Z)$ is finite, and we
put $c^{}_2 = \max_{z \in \hat l (Z)} \, |z|$.

Consider a map $\tau : X \to \R^2$ defined by
$\tau : x \mapsto ( \xi (x), \, \eta (x) )$, where
$\xi (x) = \mbox{Re} \, (l(x))$ and $\eta (x) = |\!| x |\!|$.
Let $\beta : \R^2 \setminus \{ \, 0 \, \} \to [0, 2 \pi[$
be the standard polar angle.  Note that the image of $\tau$
is contained in the sector $\{ \, \pi / 4 \leq \beta \leq
3 \pi / 4 \, \} \cup \{ 0 \}$.  Let $C \subset \R^2$ be the set of
critical values of $\tau$.  By the general topological finiteness
results of real algebraic geometry (see [GM3, Theorem 1.7]), there
exists a number $c^{}_3 > 0$ such that the set
$C \cap \{ \, (\xi, \eta) \in \R^2 \; |
\; |\!|(\xi, \eta)|\!| \geq c^{}_3 \, \}$
consists of finitely many disjoint, connected smooth curves
$\{ C_j \}_{j = 1}^{m}$, each approaching a point at
infinity in the standard projective compactification of $\R^2$
(see Figure 1).  Let $\gamma_j : \R_+ \, \to C_j$ be the length
parameterization of $C_j$, so that $\gamma_j (0)$ is the starting
point of $C_j$ (with $|\!| \gamma_j (0) |\!| = c^{}_3$), and
$|\!| \gamma'_j (t) |\!| = 1$ for all $t$.  Let $\beta_j =
\lim_{t \to \infty} \, \beta (\gamma_j (t)) \in [\pi / 4,
3 \pi / 4]$.  Renumber the curves $\{ C_j \}$ so that $\beta_j =
\pi / 2$ for $j = 1, \dots , m'$, and $\beta_j \neq \pi / 2$
for $j = m' + 1, \dots , m$.  By choosing the number $c^{}_3$
sufficiently large, we can assume that there exists a
$b > 1$ such that
$$\eta \, (\gamma_j (t)) < b \cdot |\xi (\gamma_j (t))|,$$
for any $j = m'+1, \dots , m$, and any $t \in \R_+$.

Thus, we have fixed positive constants $k, b, c^{}_1, c^{}_2,
c^{}_3$.  We now prove two technical lemmas, both direct
consequences of the estimate (\ref{tti}).

\begin{lemma}\label{tech1}
We have $| \xi (\gamma_j (t)) | < 4 k$, for any
$j = 1, \dots , m'$, and any $t \in \R_+$.
\end{lemma}
\begin{pf}
Fix $j \in \{ \, 1, \dots m' \, \}$, and write $\gamma =
\gamma_j : \R_+ \to \R^2$.  The curve $\gamma$ satisfies the
following four conditions.

(i)   $\lim_{t \to \infty} \, | \gamma (t) | = \infty$.  

(ii)  $\lim_{t \to \infty} \, \frac{\gamma (t)}{|\gamma (t)|}
      = (0, 1) \in \R^2$.

(iii) $| \gamma' (t) | = 1$, for all $t \in \R_+$.
 
(iv)  $| \, \gamma' (t) - \frac{\gamma (t)}{|\gamma (t)|} \, |
      \leq \frac{2k}{|\gamma (t)|} \,$, for all $t \in \R_+$.

\noindent
Conditions (i) - (iii) follow from the construction, and condition
(iv) follows from (\ref{tti}).

It is an elementary (but interesting) exercise in plane geometry
to verify that any curve satisfying conditions (i) - (iv) above
must also satisfy $| \xi (\gamma (t)) | < 4 k$, for all
$t \in \R_+$.
\end{pf}

\begin{lemma}\label{tech2}
Let $x \in X$ be any point with $|\!| l (x) |\!| \geq c^{}_1$.
Then there exists a tangent vector $u \in T_x X$ such that
$d \eta (u) = 0$, $d l (u) \neq 0$, and the angle
$\omega$ between the vectors $l (y)$ and $d l (u)$ in $\C$
satisfies:
$$| \, \omega - \pi / 2 \, | < \frac{1}{100} \, .$$
\end{lemma}
\begin{pf}
Let $v_0 \in T_x X$ be a vector minimizing the angle between
$v$ and $x$ in $V$, over all non-zero $v \in T_x X$.  Then, by
(\ref{tti}), we can take $u = i \, v_0$.
\end{pf}

\subsection{Proof of Proposition \ref{odlim}}

The proof of Proposition \ref{odlim} proceeds by appealing
to Figure 1 which depicts the critical locus $C$ of the map
$\tau : X \to \R^2$.  Fix a constant $c$ satisfying $c > c_i$,
for $i = 1, 2, 3$.  Let
$$D = \{ (\xi, \eta) \in \R^2 \; | \;
\xi > c, \; \eta > b \cdot \xi \}.$$

{\bf Claim:}  The homeomorphism type of the pair
$$( \{ \, x \in X \; | \; \xi (x) \leq \xi^{}_0 \, \}, \,
\{ \, \eta (x) \geq \eta^{}_0 \, \})$$
is independent of the numbers $\xi^{}_0, \eta^{}_0$, as long
as $(\xi^{}_0, \eta^{}_0) \in D$.

Proposition \ref{odlim} is an immediate consequence of this
claim.  The proof of the claim is an application of moving
the wall.  It is enough to verify the following three
conditions.

(i)   The function $\xi |_X$ has no critical points $x \in X$
with $\xi (x) > c$.

(ii)  The function $\eta |_X$ has no critical points
$x \in X$ with $\eta (x) > b \cdot c$.

(iii) The intersection $C \cap D$ is empty.

Conditions (i) and (ii) follow directly from the estimate
(\ref{tti}).  Condition (iii) follows form Lemma \ref{tech1} and
the definitions of Section 4.1.

\vspace{.1in}
\begin{center}
\leavevmode{}
\epsfbox{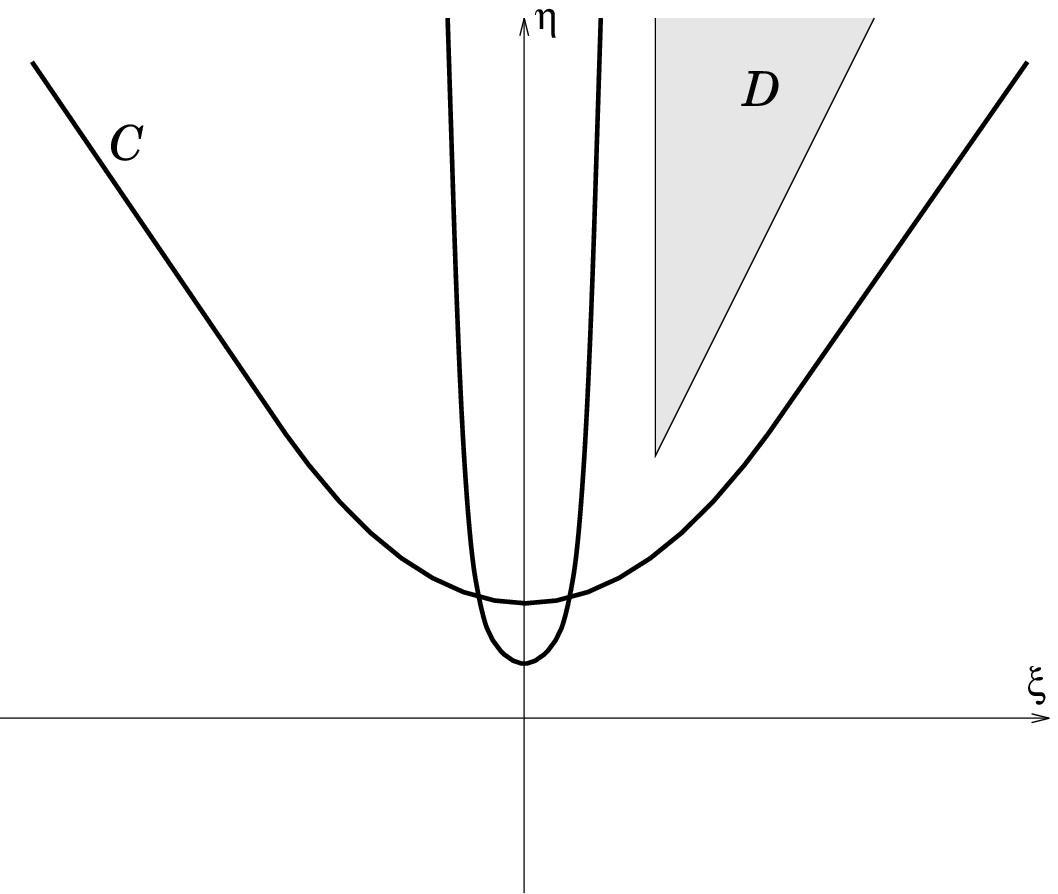}
\end{center}
\vspace{.1in}
\begin{center}
Figure 1
\end{center}
\vspace{.1in}

\subsection{Proof of Proposition \ref{auxmtw}}

Our Proof of Proposition \ref{auxmtw} will proceed in four steps.
Fix a point $(\xi^{}_0, \eta^{}_0)$ in the region $D$ of
Section 4.2.

{\bf Step 1.}  We show that there is a natural isomorphism:
$$\displaylines{\quad\quad\quad
H^* \, ( \{ \, x \in X \; | \; \xi (x) \leq \xi^{}_0
\, \}, \, \{ \, \eta (x) \geq \eta^{}_0 \, \}; \, \C ) \cong
\hfill\cr\hfill
H^* \, ( \{ \, x \in X \; | \; \xi (x) \leq \xi^{}_0
\, \}, \, \{ \, \eta (x) \geq \eta^{}_0 \, \}
\cup \{ \, \xi (x) \leq - \xi^{}_0 \, \}; \, \C ).
\quad\quad\quad}$$
By the long exact sequence of the triple, it is enough to
show that
\begin{equation}\label{HH}
H^* \, ( \{ \, x \in X \; | \; \xi (x) \leq - \xi^{}_0
\, \}, \, \{ \, \eta (x) \geq \eta^{}_0 \, \}; \, \C ) = 0.
\end{equation}
Define $X_{\leq \eta^{}_0} = \{ \, x \in X \; | \;
\eta (x) \leq \eta^{}_0 \, \}$, and $X_{\eta^{}_0} =
\{ \, x \in F \; | \; \eta (x) = \eta^{}_0 \, \}$.
The pair $(X_{\leq \eta^{}_0}, X_{\eta^{}_0})$ is a manifold
with boundary.  By excision, (\ref{HH}) is equivalent to
$$H^* \, (X_{\leq \eta^{}_0}, X_{\eta^{}_0}; \, \C) = 0.$$
This is proved as an application of Morse theory for a manifold
with boundary.  There are two conditions we have to check.

(i)   The restriction $\xi |_X$ has no critical points $x \in X$
with $\xi (x) \leq - \xi^{}_0$ and $\eta (x) \leq \eta^{}_0$.

(ii)  Let $x \in X_{\eta_0}$ be a critical point of
$\xi |_{X_{\eta_0}}$ with $\xi (x) \leq - \xi^{}_0$.
Then the differentials $d_x (\xi |_X)$ and $d_x (\eta |_X)$ in
$T^*_x X$ satisfy $d_x \xi = s \cdot d_x \eta$, for some $s < 0$.

Both (i) and (ii) follow from (\ref{tti}).

{\bf Step 2.}  Fix a smooth cut-off function $\psi : \R \to \R$
satisfying the following conditions:

(i)   $\psi (\zeta) = \xi^{}_0$, for $|\, \zeta \, |
\leq \xi^{}_0$;

(ii)  $\psi (\zeta) = 0$, for $|\, \zeta \, | \geq 5 \, \xi^{}_0$;

(iii) $0 \leq \psi (\zeta) \leq \xi^{}_0$, for all $\zeta \in \R$;

(iv)  $| \, \psi' (\zeta) \, | \leq 1/2$, for all $\zeta \in \R$.

\noindent
Let $\zeta : V \to \R$ be the imaginary part of $l$, i.e.,
$\zeta (v) = \mbox{Re} \, (- i \, l (v))$.  We show that
there is a natural isomorphism:
\begin{eqnarray}
H^* \, ( \{ \, x \in X \; | \; \xi (x) \leq \xi^{}_0
\, \}, \, \{ \, \eta (x) \geq \eta^{}_0 \, \}
\cup \{ \, \xi (x) \leq - \xi^{}_0 \, \}; \, \C ) \cong
\qquad\qquad\qquad
\nonumber \\ 
\qquad\quad
H^* \, ( \{ \, x \in X \; | \; \xi (x) \leq \psi(\zeta (x))
\, \}, \, \{ \, \eta (x) \geq \eta^{}_0 \, \}
\cup \{ \, \xi (x) \leq - \psi(\zeta (x)) \, \}; \, \C ).
\label{AA}
\end{eqnarray}
The purpose of this is to localize the problem in the $l$
directions (see Figure 2).

By the long exact sequence of the triple, it is enough to
show that
\begin{equation}
H^* \, ( \{ \, x \in X \; | \; \xi (x) \leq \xi^{}_0
\, \}, \, \{ \, \eta (x) \geq \eta^{}_0 \, \}
\cup \{ \, \xi (x) \leq - \psi(\zeta (x)) \, \}; \, \C ) = 0, \;\;
\mbox{and} \label{LL}
\end{equation}
\begin{equation}
H^* \, ( \{ \, x \in X \; | \; \xi (x) \leq - \psi(\zeta (x))
\, \}, \, \{ \, \eta (x) \geq \eta^{}_0 \, \}
\cup \{ \, \xi (x) \leq - \xi^{}_0 \, \}; \, \C ) = 0.
\label{MM}
\end{equation}
Both (\ref{LL}) and (\ref{MM}) are proved by moving the wall in
the $l$ plane.  To prove (\ref{LL}), we move the ``wall''
$\{ \, \xi (x) = \xi^{}_0 \, \}$ to
$\{ \, \xi (x) = \psi (\zeta(x)) \, \}$ through the family
$\{ \xi (x) = t \cdot \psi (\zeta(x)) + (1 - t) \cdot \xi^{}_0
\}_{\, t \in [0, 1]}$.  For $t \in [0, 1]$, define a
function $\kappa^{}_t : X \to \R$ by
$$\kappa^{}_t (x) =
\xi (x) - t \cdot \psi (\zeta(x)) - (1 - t) \cdot \xi^{}_0.$$
Then it suffices to verify the following two conditions.

(i)   Zero is not a critical value of $\kappa^{}_t$, for any
$t \in [0, 1]$.

(ii)  Zero is not a critical value of the restriction
$\kappa^{}_t |_{X_{\eta_0}}$, for any $t \in [0, 1]$
(where $X_{\eta_0} = \{ \, x \in X \; | \; \eta (x) =
\eta^{}_0 \, \}$).

Condition (i) follows from the fact that the
function $l : X \to \C$ has no critical values $z$ with
$|\, z \,| \geq \xi^{}_0$.  Condition (ii) follows from Lemma
\ref{tech2}.  Equality (\ref{MM}) is proved similarly.

\vspace{.1in}
\begin{center}
\leavevmode{}
\epsfbox{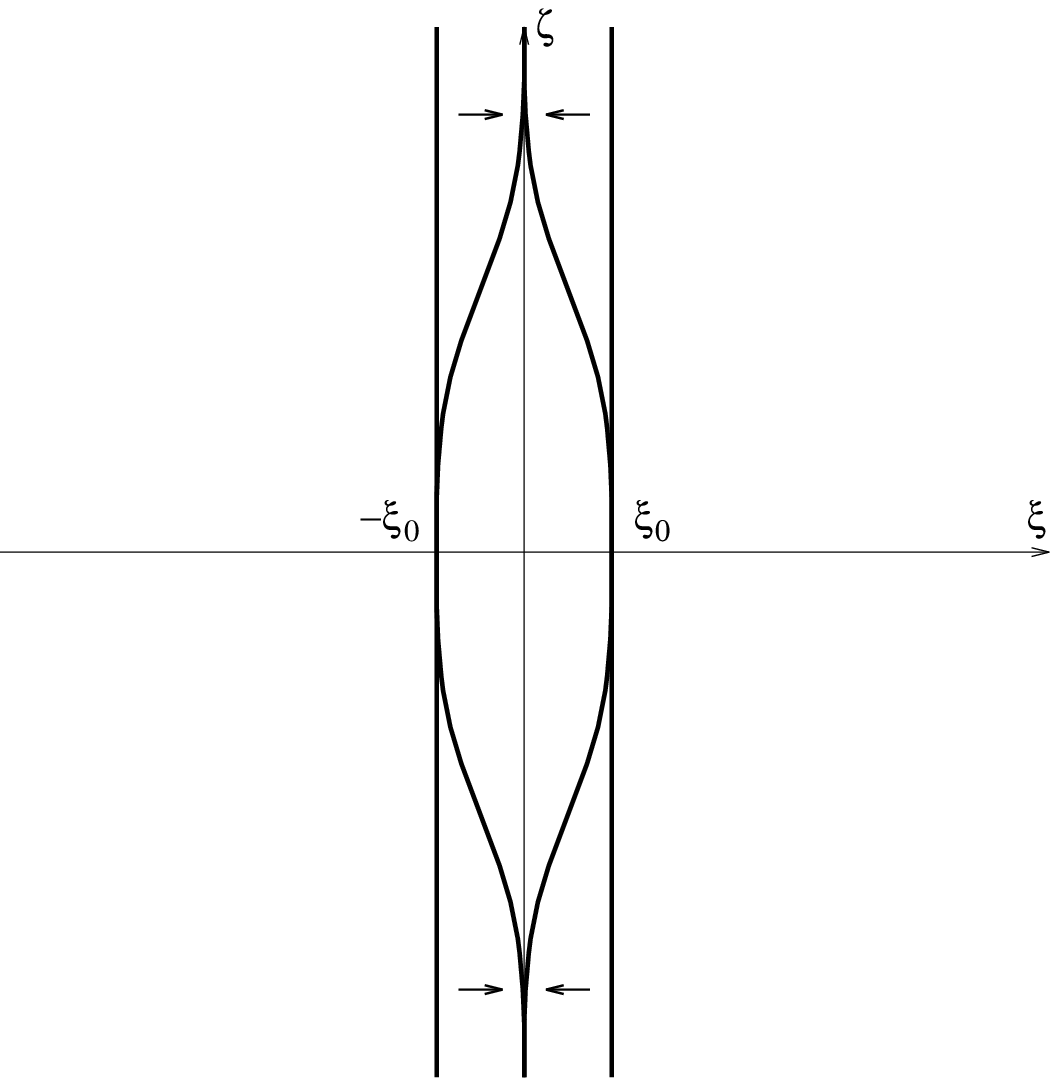}
\end{center}
\vspace{.1in}
\begin{center}
Figure 2
\end{center}
\vspace{.1in}

{\bf Step 3.}  From Steps 1 and 2, we see that the cohomology in
(\ref{AA}) is canonically independent of $\eta^{}_0$ for large
$\eta^{}_0$.  It follows that this cohomology is equal to
$$\HH^{*-n} \, ( \{ \, x \in \hat X \; | \; \bar \xi (x) \leq
\psi(\hat \zeta (x)) \, \}, \, \{ \, \hat \xi (x) \leq -
\psi(\hat \zeta (x)) \, \}; \, j^{}_! \, \C^{}_X \, [n] \, ),$$
where $\hat \xi$ and $\hat \zeta$ are the real and imaginary
parts of $\hat l$.

{\bf Step 4.}  The last step is to identify the hypercohomology
in Step 3 with the hypercohomology
$$\HH^{*-n} \, ( \{ \, x \in \hat X \, | \, \hat \xi (x) \leq
\xi^{}_0 \, \}, \, \{ \, |\!| \hat l (x) |\!| \geq 2 \xi^{}_0
\, \}; \; j^{}_! \, \C^{}_X \, [n] )$$
of Proposition \ref{auxmtw}.  This is again done by a
deformation in the $\hat l$ plane, using the fact that
$\hat l \, |_{\hat X}$ has no critical values $z$ with
$| \, z \, | \geq \xi^{}_0$ (see [GM3, Part I, Theorem 1.5]).

\section{Whitney Conditions}

We now furnish the remaining piece of the proof of Theorem
\ref{main}: a proof of Proposition \ref{whitney}. 

First, we note that the decomposition
$$\hat X^{\infty, \circ} = \bigcup_{S \in \cS^\circ \atop
K \in \cK (S)} K$$
is a Whitney stratification of $\hat X^{\infty, \circ}$.
This follows from Lemma \ref{noncrit} and the fact that
Whitney conditions are preserved under stratified transverse
intersection (see [GM3, Part I, Chapter 1.2]).  Thus, the
only non-trivial part of Proposition \ref{whitney} is the
claim that Whitney conditions (a) and (b) hold for the pair
$(X, K)$, for each $S \in \cS^\circ$ and $K \in \cK (S)$.

To check Whitney (a), consider a sequence of points
$\{ x_i \}$ in $X$ converging to $\hat x \in K \subset
\hat X$.  Let $\tilde x_i = x_i / |\!| x_i |\!|$.  By
passing if necessary to a subsequence, we may assume that
there exists a limit $q = \lim_{i \to \infty} \tilde x_i
\in V$.  Then $q$ must be a point of transverse intersection
of the stratum $S \subset \as (X)$ and the hyperplane
$\Delta \subset V$.  Whitney condition (a) for the sequence
$\{ x_i \}$ now follows from Thom's $A_f$ condition for the
sequence $\{ \tilde x_i \}$ (see assumption (ii) on the
stratification $\cS$ in Section 3.2).

It remains to check Whitney condition (b).  Fix a sequence
$\{ x_i \}$ in $X$ with a limit $\hat x \in K \subset
\hat X$.  Without loss of generality, we may assume that
$\hat l (\hat x) = 0$, so that $\hat x$ is given as
$(\Gamma, 0)$, where $\Gamma \subset \as (X) \cap \Delta$ is
a line with $\Gamma \setminus \{ 0 \} \subset S$.  We now
introduce linear coordinates $(z_1, \dots, z_d) : V \to \C^d$
such that:

(i)    the line $\Gamma$ is the $z_1$-axis;

(ii)   the covector $l = z_s$, where $s = \dim S$;

(iii)  let $q = (1, 0, \dots, 0) \in \Gamma$, then $T_q S$
is parallel to the plane $\{ z_{s+1} = z_{s+2} = \dots =
z_d = 0 \} \subset V$.

\noindent
Choose a small neighborhood $\hat U$ of $\hat x$ in $\hat V$,
and let $U = \hat U \cap V$.  The functions
$$\biggl( \frac{1}{z_1}, \; \frac{z_2}{z_1}, \; \dots, \;
\frac{z_{s-1}}{z_1}, \; z_s, \; \frac{z_{s+1}}{z_1}, \;
\dots, \; \frac{z_d}{z_1} \biggr) : \,  U \to \C^d$$
extend to give a system of coordinates on $\hat U$.
Denote these coordinates by
$$(\hat z_1, \, \dots, \, \hat z_d) : \, \hat U \to \C^d.$$
Let $\hat \pi : \hat U \to K$ be the projection defined by
$\hat z_i (\hat \pi (\hat u) ) = \hat z_i (\hat u)$, for any
$\hat u \in \hat U$ and $s \leq i \leq d$ (provided $\hat U$
is small enough, this defines $\hat \pi$ uniquely).

In the presence of Whitney condition (a), it is enough
to verify Whitney (b) for the sequences $\{ x_i \in X \}$
and $\{ \hat \pi (x_i) \in K \}$.  To do this, consider
a small neighborhood $\bar U \subset \bar V$ of the point 
$\bar x \in \bar V$ given by the line $\Gamma$.  The
functions
$$\biggl( z_1, \, \frac{z_2}{z_1}, \, \dots, \,
\frac{z_d}{z_1} \biggr) : \, \bar U \cap V \to \C^d$$
extend to give a system of coordinates $(\bar z_1, \, \dots,
\, \bar z_d)$ on $\bar U$.  Let $\bar \pi : \bar U \to
S^\infty$ be the projection defined by
$\bar z_i (\bar \pi (\bar u) ) = \bar z_i (\bar u)$, for any
$\bar u \in \bar U$ and $s \leq i \leq d$.  Then Whitney (b)
for the sequences $\{ x_i \in X \}$ and
$\{ \hat \pi (x_i) \in K \}$ follows from Whitney (b)
for the sequences $\{ x_i \in X \}$ and
$\{ \bar \pi (x_i) \in S^\infty \}$ (see assumption (iii) on
the stratification $\cS$ in Section 3.2).  This completes
the proof of Proposition \ref{whitney} and, with it, of
Theorem \ref{main}.

\section{Products of Linear Forms}

We now return to Example \ref{linear_forms} of Section 1.
Let $f  = \prod_{i \in I} l_i: V \to \C$ be a product of
linear forms $l_i \in V^*$, running over a finite set $I$,
and let $X = f^{-1} (1)$.  Then $\as (X) = f^{-1} (0)$, and
$P = \psi_f \, \C_X [d-1]$.  Singularities of $f$ have
been a subject of much recent study (see [CS] and its
references).  The following proposition shows that Theorem
\ref{main} applies to this example.  Fix a Hermitian metric
on $V$.  Let $K = \bigcap_{i \in I} \ker l_i \subset V$.

\begin{prop}\label{hyperplanes}
(i)    The norm $|\!| d_x f |\!|$ is bounded from below
on $X$.

(ii)   The variety $X$ is transverse to infinity.

(iii)  The variety $Y = \overline{p_2 (T_X^* V)}$, associated
to $X$ as in Theorem \ref{main}, is the orthogonal complement
$K^\perp \subset V^*$.
\end{prop}
\begin{pf}
We argue by induction on $d$ and the cardinality $|I|$.
Case $d = 1$ or $|I| = 1$ is trivial.  Assume now $d > 1$
and $|I| > 1$.  If $K \neq \{ 0 \}$, we may reduce the
parameter $d$ by passing to the quotient $V / K$.  Assume
now $K = \{ 0 \}$.

By a standard compactness argument, parts (i) and (ii) of the
proposition will follow if we show that for any sequence
$\{ x_j \}$ in $X$, converging to a limit $\bar x \in \bar X
\subset \bar V$:

(a) the norm $|\!| d_{x_j} f |\!|$ is bounded away from zero;

(b) the product $\angle (x_j, T_{x_j} X) \cdot
|\!| x_j |\!|$ is bounded from above.

\noindent
Let $\Gamma \subset V$ be the line corresponding to $\bar x$.
Define $I' = \{ i \in I \, | \, \ker l_i \supset \Gamma \}$.
Let $m = | I \setminus I' |$; note that $m > 0$.  Set
$g = \prod_{i \in I'} l_i$, and
$h = \prod_{i \in I \setminus I'} l_i$, so that
$f = g \cdot h$.  Then we have the asymptotics (up to a
multiplicative constant):
$$|h(x_j)| \sim |\!| x_j |\!|^m, \;\; |g(x_j)| \sim
|\!| x_j |\!|^{-m}, \;\; |\!| d_{x_j} h |\!| \sim
|\!| x_j |\!|^{m-1}.$$
By the Leibniz's rule, $d f = g \cdot d h + h \cdot d g$.  By
the above asymptotics,
$$|\!| g (x_j) \cdot d_{x_j} h |\!| \sim |\!| x_j |\!|^{-1}.$$
On the other hand, by part (i) of the induction hypothesis on
$g$, the norm $|\!| h (x_j) \cdot d_{x_j} g |\!|$ is bounded
away from zero.  This immediately implies (a).
Furthermore, it shows that there is a constant  $c > 0$, such
that the angle $\alpha_j$ between the hyperplanes
$\ker d_{x_j} f = T_{x_j} X$ and $\ker d_{x_j} g$ in
$T_{x_j} V \cong V$ satisfies $\alpha_j < c / |\!| x_j |\!|$,
for all $j$.  Together with part (ii) of the induction
hypothesis on $g$, this estimate implies (b).

To prove part (iii) of the proposition, denote by $\cV$ the
linear stratification of $V$ underlying the arrangement
$\as (X)$, and by $\cV^*$ the dual stratification of $V^*$.
Note that $\cV^*$ has more than one codimension one stratum.
An inductive argument, as in the proof of parts (i) and (ii),
shows that $Y$ contains every codimension one stratum of
$\cV^*$.  But $Y$ is irreducible, therefore it must be all
of $V^*$.
\end{pf}

We conclude with an observation which gives, in this example,
an algorithm for computing the rank of the local system $\cL$
of Theorem \ref{main} from the combinatorics of the
arrangement $\as (X)$.  Let $\cV$ be as in the proof of
Proposition \ref{hyperplanes}, and $\cS = \{ S \in \cV \, |
\, S \subset \as (X) \}$ be the induced stratification of
$\as (X)$.  Note that $K$ is the smallest stratum of $\cS$.
For $S \in \cS$, let $m (S)$ be the multiplicity of the
conormal bundle $T^*_S V$ in the characteristic cycle of $P$.
Also, let $c (S)$ be the (complex) codimension of $S$ in
$\as (X)$.  Note that $\rank \cL = m (K)$.

\begin{prop}\label{rankl}
Let $\chi (X)$ be the Euler characteristic of $X$.  Then
$$\chi (X) = \sum_{S \in \cS} (-1)^{c(S)} \cdot m (S).$$
\end{prop}
\begin{pf}
This follows from the Morse theory of the distance to a
generic point in $V$, as in [GM3, Part III].
\end{pf}

Proposition \ref{rankl} may be used to compute the numbers
$m (S)$ inductively.  The only ingredient for the induction
step is the Euler characteristic $\chi (X)$.  To compute it,
note that $X$ is an $|I|$-fold cover over
$(V \setminus \as (X)) / \C^*$ (see [CS]).  Therefore,
$\chi (X) = |I| \cdot \chi((V \setminus \as (X)) / \C^*)$.
Pick any $i \in I$, and let
$H_i = \{ v \in V \, | \, l_i (v) = 1 \}$.  Then
$$(V \setminus \as (X)) / \C^* \cong
H_i \setminus (\as (X) \cap H_i).$$
Finally, the Euler characteristic (in fact, the full
cohomology) of the complement $H_i \setminus (\as (X) \cap
H_i)$ may be computed as in [GM3, Part III].

\vspace{.1in}

\noindent
Department of Mathematics, MIT, 77 Massachusetts Ave.,
Cambridge, MA 02139

\noindent
{\it grinberg@math.mit.edu}

\end{document}